\documentclass[10pt,twoside]{article}

\usepackage{amsmath,amssymb,amsthm,hyperref,graphicx,color}

\numberwithin{equation}{section}
\numberwithin{table}{section}
\numberwithin{figure}{section}

\theoremstyle{plain}
\newtheorem{theorem}{Theorem}[section]

\newtheorem{corollary}[theorem]{Corollary}
\newtheorem{remark}[theorem]{Remark}

\def\l|{\left\vert}
\def\r|{\right\vert}

\def\lef|{\left\vert}
\def\llef|{\Bigl\vert}

\def\({\left(}
\def\){\right)}
\def\[{\left[}
\def\]{\right]}
\def\ll|{\Bigl\vert}
\def\rr|{\Bigr\vert}

\def\dfrac#1#2{{\displaystyle\frac{#1}{#2}}}

\def\toD{\stackrel{\cal{L}}{\to}}

\newcommand\eqdist{\stackrel{\mathcal{L}}{=}}
\newcommand\wtilde{\widetilde}

\def\bmu{{\boldsymbol\mu}}

\def\bmuhat{{\widehat{\boldsymbol\mu}}}

\def\bLambda{{\boldsymbol\Lambda}}

\def\bSigma{{\boldsymbol\Sigma}}
\def\bSigmahat{{\widehat{\boldsymbol\Sigma}}}
\def\bDeltahat{{\widehat{\boldsymbol\Delta}}}

\def\bnu{{\boldsymbol\nu}}

\def\bartau{\bar{\tau}}

\def\tr{{\rm{tr}}}
\def\Cov{{\rm{Cov}}}
\def\Var{{\rm{Var}}}

\def\bzero{\boldsymbol{0}}

\def\b1{\,\boldsymbol{1}}

\def\bA{{\boldsymbol{A}}}
\def\bB{{\boldsymbol{B}}}
\def\bC{{\boldsymbol{C}}}

\def\bI{{\boldsymbol{I}}}

\def\bS{{\boldsymbol{S}}}

\def\bX{{\boldsymbol{X}}}
\def\bXhat{\widehat{\boldsymbol{X}}}
\def\bY{{\boldsymbol{Y}}}

\def\bW{{\boldsymbol{W}}}

\def\bZ{{\boldsymbol{Z}}}
\def\bXi{{\boldsymbol{\Xi}}}
\def\bTheta{{\boldsymbol{\Theta}}}
\def\bomega{{\boldsymbol{\omega}}}

\begin{document}

\title{\bf\large Kurtosis Tests for Multivariate Normality with Monotone Incomplete Data}

\author{{Tomoya Yamada,}\thanks{Department of Economics, Sapporo Gakuin University, 11 Bunkyodai, Ebetsu, Hokkaido, Japan.  Research supported by a sabbatical leave-of-absence from 
Sapporo Gakuin University.  E-mail address: tomoya@sgu.ac.jp}
\ \ {Megan M. Romer,}\thanks{Department of Statistics, Penn State University, University Park, PA 16802, U.S.A.  Research supported in part by National Science Foundation grant DMS-1309808.  E-mail address: mmd198@stat.psu.edu}
\ \ {and \ Donald St. P. Richards}\thanks{Department of Statistics, Penn State University, University Park, PA 16802, U.S.A.  Research supported in part by National Science Foundation grant DMS-1309808.  E-mail address: richards@stat.psu.edu
\endgraf
\endgraf
}
}

\maketitle

\begin{abstract}
\footnotesize{
We consider the problem of testing multivariate normality when the data consists of a random sample of two-step monotone incomplete observations.  We define for such data a generalization of Mardia's statistic for measuring kurtosis, derive the asymptotic non-null distribution of the statistic under certain regularity conditions and against a broad class of alternatives, and give an application to a well-known data set on cholesterol measurements.
}
\end{abstract}

\section{Introduction}
\label{introduction}
\setcounter{equation}{0}

In many statistical applications, and especially in epidemiology and biostatistics, incomplete data arise for a variety of reasons; cf., Eaton and Kariya (1983), Garren and Peddada (2000), Little and Rubin (2002), Peddada, Harris, and Davidov (2010), Krishnamoorthy and Yu (2012), and Davidov and Peddada (2013).  Consequently, much work has been done on explicit formulas that allow for statistical inference with incomplete data to achieve specified levels of significance.  

In this paper, we consider the problem of kurtosis tests for multivariate normality using two-step monotone incomplete data.  Following our earlier work on monotone incomplete multivariate normal data (Chang and Richards, 2009, 2010; Richards and Yamada, 2010; Romer, 2009; Romer and Richards, 2010, 2013; Yamada, 2013), we write the data in the form 
\begin{equation}
\label{monotonesample}
\binom{\bX_1}{\bY_1} \ 
\binom{\bX_2}{\bY_2} \ 
\cdots \ 
\binom{\bX_n}{\bY_n} \ \; 
\begin{matrix}        \\ \bY_{n+1} \end{matrix} \ \
\begin{matrix}        \\ \bY_{n+2} \end{matrix} \ \
\begin{matrix}        \\ \cdots \end{matrix}  \ \
\begin{matrix}        \\ \bY_N \end{matrix}
\end{equation}
where each $\bX_j$, $1 \le j \le n$, is $p \times 1$ and each $\bY_j$, $1 \le j \le N$, is $q \times 1$.  

As in our earlier work, we assume that the data are missing completely at random (MCAR).  It follows from results of Eaton and Kariya (1983), Hao and Krishnamoorthy (2001), and others that an explicit solution for the likelihood equations for the covariance matrix requires the MCAR assumption.  Thus, the MCAR assumption and the monotone data pattern ensure the validity of likelihood inference and a unique, explicit solution to the likelihood equations.  

In this paper, we define for data of the form (\ref{monotonesample}) a generalization of Mardia's statistic for testing kurtosis.  We derive the asymptotic non-null and null distributions of the new statistic under certain regularity conditions.  We apply our results to a well-known Pennsylvania cholesterol data set (Ryan, Joiner, and Cryer, 2005, p.~267) which has been used widely to illustrate statistical methods for analyzing monotone incomplete multivariate data.  Our results provide an invariant test of normality for that data, adding to the literature on that subject (see Henze (2002) for an extensive account of invariant testing procedures in the complete case), and we complement results of Romer (2009) and Romer and Richards (2013), where exploratory methods were applied to the cholesterol data set.  

Our results are as follows.  In Section \ref{prelim}, we provide the background for testing kurtosis with the data (\ref{monotonesample}).  We define in Section \ref{kurtosis} the kurtosis statistic, $b_{2,p,q}$, and prove that $b_{2,p,q}$ is identical to a statistic constructed from the observed data and from data imputed by linear regression methods.  Further, we establish an invariance property of $b_{2,p,q}$ which allows us to reduce the general problem to the canonical case in which the population has mean zero and identity covariance matrix.  We state in Section \ref{asymptoticdistns} the null and non-null asymptotic distributions of $b_{2,p,q}$ corresponding, respectively to the normal case and to a broad class of alternatives defined by moment conditions on the distribution, and we apply those results to the Pennsylvania cholesterol data, reaching conclusions similar to Romer (2009).  We derive in Appendix \ref{Sigma-asymptotics} asymptotic expansions of $\bSigmahat$ and $\bSigmahat^{-1}$, where $\bSigmahat$ is an estimator of $\bSigma$, the covariance matrix of the population underlying the sample (\ref{monotonesample}), and then in Appendix \ref{b1b2asymptotics} we apply those expansions to derive asymptotic expansions of two statistics used to construct $b_{2,p,q}$.  Finally, we obtain in Appendix \ref{basymptotics} the null and non-null asymptotic distributions of $b_{2,p,q}$.

\section{Notation and preliminaries}
\label{prelim}
\setcounter{equation}{0}

Throughout the paper, we follow the notation of Chang and Richards (2009, 2010).  Thus, all matrices and vectors are written in boldface type; $\bI_d$ denotes the identity matrix of order $d$; and $\bzero$ denotes any matrix or vector of zeros, the dimension of which will be clear from the context.  We also let $\tau = n/N$ denote the proportion of observations in (\ref{monotonesample}) that are complete, and set $\bartau = 1-\tau \equiv (N-n)/N$.  

Define the sample means 
\begin{equation}
\label{samplemeans}
\bar \bX = \frac{1}{n} \sum_{j=1}^n \bX_j, \qquad 
\bar \bY_1 = \frac{1}{n} \sum_{j=1}^n \bY_j, \qquad 
\bar \bY_2 = \frac{1}{N-n}\sum_{j=n+1}^N \bY_j,
\end{equation}
and 
\begin{equation}
\label{grandmean}
\bar \bY = \frac1N \sum_{j=1}^N \bY_j \equiv \tau \bar \bY_1+\bartau \bar \bY_2.
\end{equation}
We also define corresponding matrices of sums of squares and products, 
\begin{equation}
\label{sumsqmatrices}
\begin{array}{ll}
\bA_{11} = \sum\limits_{j=1}^n (\bX_j-\bar \bX)(\bX_j-\bar \bX)', &
\bA_{12} = \bA_{21}' = 
\sum\limits_{j=1}^n (\bX_j-\bar \bX)(\bY_j-\bar \bY_1)', \\
\bA_{22,n} = \sum\limits_{j=1}^n (\bY_j-\bar \bY_1)(\bY_j-\bar \bY_1)', &
\bA_{22,N} = \sum\limits_{j=1}^N (\bY_j-\bar \bY)(\bY_j-\bar \bY)',
\end{array}
\end{equation}
set $\bA_{11\cdot2,n}=\bA_{11}-\bA_{12}\bA_{22,n}^{-1}\bA_{21}$, and let 
\begin{equation}
\label{Amatrix}
\bA = \begin{pmatrix}\bA_{11} &\bA_{12}\\ \bA_{21}&\bA_{22,n}\end{pmatrix}.
\end{equation}

Let $\bmu$ denote the mean and $\bSigma$ the covariance matrix of the population underlying the data (\ref{monotonesample}); we assume that $\bSigma$ is nonsingular.  We partition $\bmu$ and $\bSigma$ similar to (\ref{monotonesample}), so that 
$$
\bmu = \begin{pmatrix} \bmu_1 \\ \bmu_2 \end{pmatrix}, 
\qquad
\bSigma = \begin{pmatrix}
\bSigma_{11} & \bSigma_{12} \\ \bSigma_{21} & \bSigma_{22}
\end{pmatrix},
$$
where $\bmu_1$ and $\bmu_2$ are $p \times 1$ and $q \times 1$ vectors, respectively, and the submatrices $\bSigma_{11}$, $\bSigma_{12} = \bSigma_{21}'$, and $\bSigma_{22}$ are of order 
$p \times p$, $p \times q$, and $q \times q$, respectively.  We also define 
$$
\bmuhat = \begin{pmatrix}\bmuhat_1 \\ \bmuhat_2\end{pmatrix}, 
\qquad
\bSigmahat = \begin{pmatrix} 
\bSigmahat_{11} &\bSigmahat_{12} \\ 
\bSigmahat_{21}&\bSigmahat_{22} 
\end{pmatrix}, 
$$
where 
\begin{equation}
\label{muhat}
\bmuhat_1 = \bar \bX - 
\bartau \bA_{12}\bA_{22,n}^{-1}(\bar \bY_1-\bar \bY_2),
\qquad 
\bmuhat_2 = \bar{\bY},
\end{equation}
and 
\begin{equation}
\begin{aligned}
\label{Sigmahat}
& \bSigmahat_{11} = \frac1n \bA_{11\cdot2,n} 
+ \frac1N \bA_{12}\bA_{22,n}^{-1}\bA_{22,N}\bA_{22,n}^{-1}\bA_{21}, \\
& \bSigmahat_{12} = \bSigmahat_{21}' 
= \frac1N \bA_{12}\bA_{22,n}^{-1}\bA_{22,N}, 
\quad 
\bSigmahat_{22} = \frac1N \bA_{22,N}.
\end{aligned}
\end{equation}

If the underlying population is multivariate normal, denoted by $N_{p+q}(\bmu,\bSigma)$, then $\bmuhat$ and $\bSigmahat$ are the maximum likelihood estimators of $\bmu$ and $\bSigma$, respectively.  We refer to Chang and Richards (2009, 2010), Richards and Yamada (2010), Romer (2009), and Romer and Richards (2010) for results on the distributions of $\bmuhat$ and $\bSigmahat$, and inference for $\bmu$ and $\bSigma$.  

As noted by Chang and Richards (2009, p. 1886), the statistical model underlying two-step monotone incomplete multivariate normal data is related to double sampling designs, in which additional data are collected on a subset of variables in order to improve estimation of a parameter; cf. Little (1976, p. 594) and Cohn, Davidov, and Haitovsky (2008).

\section{Testing kurtosis with monotone incomplete data}
\label{kurtosis}
\setcounter{equation}{0}

\subsection{The kurtosis statistic}
\label{kurtosisstatistic}

Consider a $d$-dimensional multivariate population represented by a random vector $\bZ$ with mean vector $\bmu$ and nonsingular covariance matrix $\bSigma$, and kurtosis parameter 
$$
\beta_{2,d} = E[(\bZ-\bmu)'\bSigma^{-1}(\bZ-\bmu)]^2.
$$
To perform inference for $\beta_{2,d}$ with the monotone incomplete data (\ref{monotonesample}), we define the statistic 
\begin{equation}
\begin{aligned}
\label{b2pq}
b_{2,p,q} = \frac1N\Bigg\{c_1 \sum_{j=1}^n \Bigg[ & 
\Bigg(\binom{\bX_j}{\bY_j}-\bmuhat\Bigg)'\bSigmahat^{-1}
\Bigg(\binom{\bX_j}{\bY_j}-\bmuhat\Bigg)\Bigg]^2 \\ 
& + \, c_2 \sum_{j=n+1}^N \big[
(\bY_j-\bmuhat_2)'\bSigmahat^{-1}_{22}
(\bY_j-\bmuhat_2)\big]^2\Bigg\},
\end{aligned}
\end{equation}
where $c_1, c_2 > 0$ are constants.  In general, $c_1$ and $c_2$ can depend on $n$ and $N$ subject to the conditions $c_1 = O(\tau)$, $c_2 = O(\bartau)$, and $c_1, c_2 \not\to 0$ as $n, N \to \infty$; e.g., $(c_1,c_2) = (\tau,\bartau)$ with $n/N \to \delta \in (0,1)$.  Alternatively, $(c_1,c_2)$ can be chosen to minimize $\sigma^2$, the asymptotic variance of $b_{2,p,q}$ under the null hypothesis; this can be obtained by minimizing, over all $(c_1,c_2)$ subject to a suitable constraint on $c_1$ and $c_2$, a formula for $\sigma^2$ in (\ref{nullsigmasq}).  

Each term in (\ref{b2pq}) is an analog of the well-known statistic of Mardia (1970, 1974) for testing kurtosis with complete data, and our usage of different weights is due to the fact that the incomplete data $\bY_j$, $j=n+1,\ldots,N$ provide partial information about the population.  

We may also motivate the statistic $b_{2,p,q}$ as follows: First, we impute each missing observation $\bX_j$, $j=n+1,\ldots,N$, using a linear regression imputation scheme, 
\begin{equation}
\label{Xhatj}
\bXhat_j = \widehat{E}(\bX_j|\bY_j) \equiv 
\bmuhat_1 + \bSigmahat_{12} \bSigmahat_{22}^{-1} 
({\bY}_{j}-\bmuhat_2),
\end{equation}
which is motivated by the formula for the conditional expectation of a partitioned multivariate normally distributed random vector.  Under the hypothesis of multivariate normality, $\bXhat_j$ is the maximum likelihood estimator of $E(\bX_j|\bY_j)$, the conditional expectation of $\bX_j$ given $\bY_j$.  Second, we use as our data the merged sets of observed and imputed data vectors, 
\begin{equation}
\label{impute}
\binom{\bX_1}{\bY_1} \ 
\binom{\bX_2}{\bY_2} \ 
\cdots \ 
\binom{\bX_n}{\bY_n} \ 
\binom{\,\bXhat_{n+1}}{\bY_{n+1}} \ \
\binom{\,\bXhat_{n+2}}{\bY_{n+2}} \ \
\cdots \ 
\binom{\,\bXhat_N}{\bY_N}.
\end{equation}
To perform inference about $\beta_{2,d}$, it is natural to use the statistic
\begin{equation}
\begin{aligned}
\label{bhat}
\widehat b_{2,p,q} = \frac1N\Bigg\{& c_1 \sum_{j=1}^n \Bigg[
\Bigg(\binom{\bX_j}{\bY_j}-\bmuhat\Bigg)'\bSigmahat^{-1}
\Bigg(\binom{\bX_j}{\bY_j}-\bmuhat\Bigg)\Bigg]^2 \\ 
& + \ c_2 \sum_{j=n+1}^N \Bigg[
\Bigg(\binom{\;\bXhat_j}{\bY_j}-\bmuhat\Bigg)'\bSigmahat^{-1}
\Bigg(\binom{\;\bXhat_j}{\bY_j}-\bmuhat\Bigg)\Bigg]^2\Bigg\},
\end{aligned}
\end{equation}
an analog of Mardia's statistic based on the vectors in (\ref{impute}).  We again use possibly different weights, $c_1$ and $c_2$, to reflect the fact that some data are imputed, hence they provide less information about $\bmu$ and $\bSigma$ than in the case in which all observations are fully observed.  It is remarkable that $b_{2,p,q} \equiv \widehat{b}_{2,p,q}$, a result established in Theorem \ref{invariant} given in Section \ref{invariance}.

\vskip5pt

Set $\bZ_j = \displaystyle\binom{\bX_j}{\bY_j}$, $j=1,\ldots,n$; then, (\ref{b2pq}) becomes 
$$
b_{2,p,q} = \frac1N\Big\{c_1 \sum_{j=1}^n \big[(\bZ_j-\bmuhat)'\bSigmahat^{-1}(\bZ_j-\bmuhat)\big]^2 + c_2 \sum_{j=n+1}^N \big[(\bY_j-\bmuhat_2)'\bSigmahat^{-1}_{22}(\bY_j-\bmuhat_2)\big]^2\Big\}.
$$
Also, let 
\begin{equation}
\label{z1bar}
\bar\bZ_1 = \frac1n\sum_{j=1}^n \bZ_j \equiv \binom{\bar\bX}{\bar\bY_1},
\end{equation}
and 
\begin{equation}
\label{ytilde}
\widetilde \bY = 
\begin{pmatrix} 
\bA_{12}\bA_{22,n}^{-1}(\bar\bY_1 - \bar\bY_2) \\ 
\bar\bY_1 - \bar\bY_2
\end{pmatrix}
\equiv \bA 
\begin{pmatrix} \bzero & \bzero \\ 
\bzero & \bA_{22,n}^{-1} \end{pmatrix}
\begin{pmatrix} \bzero \\ \bar \bY_1-\bar \bY_2 \end{pmatrix}.
\end{equation}
By a direct calculation using (\ref{Amatrix}) and (\ref{muhat}), we deduce that 
\begin{equation}
\label{zjminusmuhat}
\bZ_j - \bmuhat = \bZ_j - \bar\bZ_1 + \bartau \widetilde\bY.
\end{equation}

\subsection{An invariance property of \texorpdfstring{$\boldsymbol{b_{2,p,q}}$}{b2pq}}
\label{invariance}

Define the statistics 
\begin{align*}
b_{2,p,q}^{(1)} & = \sum_{j=1}^n \big[(\bZ_j-\bmuhat)'\bSigmahat^{-1}(\bZ_j-\bmuhat)\big]^2, \\
\intertext{and} 
b_{2,p,q}^{(2)} & = \sum_{j=n+1}^N\big[(\bY_j - \bmuhat_2)'\bSigmahat^{-1}_{22}(\bY_j-\bmuhat_2)\big]^2.
\end{align*}
Then, 
\begin{equation}
\label{b-formula}
b_{2,p,q} = \frac1N\Big(c_1b_{2,p,q}^{(1)} + c_2b_{2,p,q}^{(2)}\Big).
\end{equation}

Let $\bLambda_{11}$ and $\bLambda_{22}$ be $p \times p$ and $q \times q$ positive definite matrices, respectively; let $\bLambda_{12}$ be a $p \times q$ matrix; and let $\bnu_1$ and $\bnu_2$ be $p \times 1$ and $q \times 1$ vectors, respectively.  Set 
$$
\bLambda = 
\begin{pmatrix} 
\bLambda_{11} & \bzero \\ \bzero & \bLambda_{22} 
\end{pmatrix},
\quad 
\bC = 
\begin{pmatrix} 
\bI_p & \bLambda_{12} \\ \bzero & \bI_q 
\end{pmatrix},
\quad
\bnu = \begin{pmatrix}
\bnu_1 \\\bnu_2
\end{pmatrix},
$$
and consider the group of affine transformations of the data (\ref{monotonesample}) of the form 
\begin{equation}
\label{affine}
\begin{array}{ccll}
\displaystyle{\binom{\bX_j}{\bY_j}} & \to & \bLambda \bC \displaystyle{\binom{\bX_j}{\bY_j}} + \bnu, 
& j=1,\ldots,n \\
& &\\
\bY_j & \to & \bLambda_{22}\bY_j + \bnu_2, & j=n+1,\ldots,N
\end{array}
\end{equation}
Now we have the following result: 

\begin{theorem}
\label{invariant}
The kurtosis statistics $b_{2,p,q}^{(1)}$, $b_{2,p,q}^{(2)}$, and $b_{2,p,q}$ are invariant under the transformations (\ref{affine}).  Moreover, $b_{2,p,q} \equiv \widehat{b}_{2,p,q}$.
\end{theorem}

\noindent{\it Proof}.  
Under the transformation (\ref{affine}), we verify using (\ref{samplemeans}) and (\ref{grandmean}) that $\bar\bX$, $\bar\bY_1$, $\bar\bY_2$, and $\bar\bY$ are transformed to $\bLambda_{11}(\bar\bX_1 + \bLambda_{12}\bar\bY_1) + \bnu_1$, $\bLambda_{22}\bar\bY_1 + \bnu_2$, $\bLambda_{22}\bar\bY_2 + \bnu_2$, and $\bLambda_{22}\bar\bY + \bnu_2$, respectively.  Further, by (\ref{sumsqmatrices}), the matrix $\bA$ in (\ref{Amatrix}) is transformed to $\bLambda\bC\bA\bC'\bLambda$, i.e., 
\begin{eqnarray}
\label{Atransform}
\bA_{11} & \to & \bLambda_{11}(\bA_{11} + \bLambda_{12}\bA_{21} + \bA_{12}\bLambda_{21} + \bLambda_{12}\bA_{22,n}\bLambda_{21})\bLambda_{11} \nonumber \\
\bA_{12} & \to & \bLambda_{11}(\bA_{12} + \bLambda_{12}\bA_{22,n})\bLambda_{22}, \\
\bA_{22,n} & \to & \bLambda_{22}\bA_{22,n}\bLambda_{22} \nonumber
\end{eqnarray}
Hence $\bA_{11\cdot 2,n} \to \bLambda_{11}\bA_{11\cdot 2,n}\bLambda_{11}$ and $\bA_{22,N} \to \bLambda_{22}\bA_{22,N}\bLambda_{22}$.  Further, it follows from (\ref{Sigmahat}) and (\ref{Atransform}) that $\bSigmahat_{11\cdot 2}$ and $\bSigmahat_{22}$ are transformed to $\bLambda_{11}\bSigmahat_{11\cdot 2}\bLambda_{11}$ and $\bLambda_{22}\bSigmahat_{22}\bLambda_{22}$, respectively.  

By a well-known quadratic identity (Anderson, 2003, p.~63, Exercise~2.54), for $j=1,\ldots,n$, 
\begin{align}
\label{quadidentity}
(\bZ_j-\bmuhat)'&\bSigmahat^{-1}(\bZ_j-\bmuhat) \nonumber \\
\equiv \ & \begin{pmatrix} \bX_j-\bmuhat_1 \\ \bY_j-\bmuhat_2 \end{pmatrix}'
\begin{pmatrix} 
\bSigmahat_{11} &\bSigmahat_{12} \\ 
\bSigmahat_{21}&\bSigmahat_{22} 
\end{pmatrix}^{-1}
\begin{pmatrix} \bX_j-\bmuhat_1 \\ \bY_j-\bmuhat_2 \end{pmatrix} \nonumber \\
= \ & \big(\bX_j-\bmuhat_1-\bSigmahat_{12}\bSigmahat_{22}^{-1}(\bY_j-\bmuhat_2)\big)'\bSigmahat_{11\cdot 2}^{-1}\big(\bX_j-\bmuhat_1-\bSigmahat_{12}\bSigmahat_{22}^{-1}(\bY_j-\bmuhat_2)\big) \nonumber \\
& \qquad + (\bY_j-\bmuhat_2)'\bSigmahat_{22}^{-1}(\bY_j-\bmuhat_2).
\end{align}
It follows from (\ref{muhat}) that $\bX_j-\bmuhat_1-\bSigmahat_{12}\bSigmahat_{22}^{-1}(\bY_j-\bmuhat_2) = \bX_j-\bar\bX-\bA_{12}\bA_{22,n}^{-1}(\bY_j-\bar\bY_1)$, and the latter expression is transformed by (\ref{affine}) to 
\begin{align*}
\bLambda_{11}(\bX_j + \bLambda_{12}\bY_j) - \bLambda_{11}(\bar\bX + \bLambda_{12}\bar\bY_1) - \bLambda_{11}&(\bA_{12} + \bLambda_{12}\bA_{22,n})\bA_{22,n}^{-1}(\bY_j-\bar\bY_1) \\ 
& = \bLambda_{11}\big(\bX_j-\bar\bX - \bA_{12}\bA_{22,n}^{-1}(\bY_j-\bar\bY_1)\big) \\
& \equiv \bLambda_{11}\big(\bX_j-\bmuhat_1-\bSigmahat_{12}\bSigmahat_{22}^{-1}(\bY_j-\bmuhat_2)\big).
\end{align*}
Recall that $\bSigmahat_{11\cdot 2}$ is transformed to $\bLambda_{11}\bSigmahat_{11\cdot 2}\bLambda_{11}$, so it follows that the first term in (\ref{quadidentity}) remains invariant under (\ref{affine}).  

In similar fashion, (\ref{affine}) transforms $\bY_j-\bmuhat_2$ to $\bLambda_{22}(\bY_j-\bmuhat_2)$, $j=1,\ldots,N$ and, as noted earlier, also transforms $\bSigmahat_{22}$ to $\bLambda_{22}\bSigmahat_{22}\bLambda_{22}$; consequently, the second term in (\ref{quadidentity}) remains invariant under (\ref{affine}).  It follows that $b_{2,p,q}^{(1)}$ and $b_{2,p,q}^{(2)}$ each are invariant under (\ref{affine}), so is $b_{2,p,q}$.  

Finally, to show that $b_{2,p,q} \equiv \widehat{b}_{2,p,q}$, we apply the earlier quadratic identity (Anderson, 2003, loc. cit.) for $j=n+1,\ldots,N$ to obtain 
\begin{align*}
\begin{pmatrix} \bXhat_j-\bmuhat_1 \\ \bY_j-\bmuhat_2\end{pmatrix}' &
\begin{pmatrix} 
\bSigmahat_{11} &\bSigmahat_{12} \\ 
\bSigmahat_{21} &\bSigmahat_{22} 
\end{pmatrix}^{-1}
\begin{pmatrix} \bXhat_j-\bmuhat_1 \\ \bY_j-\bmuhat_2\end{pmatrix} \\
= & \ \big(\bXhat_j-\bmuhat_1-\bSigmahat_{12} \bSigmahat_{22}^{-1}
(\bY_j-\bmuhat_2)\big)'\bSigmahat_{11\cdot 2}^{-1}
\big(\bXhat_j-\bmuhat_1-\bSigmahat_{12} \bSigmahat_{22}^{-1}
(\bY_j-\bmuhat_2)\big) \\
& \ + (\bY_j-\bmuhat_2)'\bSigmahat_{22}^{-1}(\bY_j-\bmuhat_2).
\end{align*}
By the definition in (\ref{Xhatj}) of $\bXhat_j$, it follows that $\bXhat_j - \bmuhat_1 - \bSigmahat_{12}\bSigmahat_{22}^{-1}(\bY_j-\bmuhat_2) \equiv \bzero$.  Therefore the last terms in (\ref{b2pq}) and (\ref{bhat}) are identical, so we obtain $b_{2,p,q} \equiv \widehat b_{2,p,q}$.
$\qed$

\smallskip

\begin{remark}
\label{normalcase}
{\rm In the multivariate normal case, Romer and Richards (2010, Proposition 2.1) showed that the invariance of $b_{2,p,q}$ is due to the fact that $\bmuhat$ and $\bSigmahat$, being maximum likelihood estimators, are equivariant.  Specifically, $\bmuhat$ and $\bSigmahat$ are transformed under (\ref{affine}) to $\bLambda\bC\bmuhat + \bnu$ and $\bLambda\bC\bSigmahat\bC'\bLambda$, respectively; therefore the quadratic form $(\bZ_j-\bmuhat)'\bSigmahat^{-1}(\bZ_j-\bmuhat)$, $j=1,\ldots,n$ is transformed to 
\begin{align*}
\big((\bLambda\bC\bZ_j+\bnu)-(\bLambda\bC\bmuhat+\bnu)\big)'
(\bLambda\bC\bSigmahat\bC'\bLambda)^{-1}
\big((\bLambda\bC\bZ_j+\bnu)&-(\bLambda\bC\bmuhat+\bnu)\big) \\
& \equiv (\bZ_j-\bmuhat)'\bSigmahat^{-1}(\bZ_j-\bmuhat),
\end{align*}
and this proves that $b_{2,p,q}^{(1)}$ is invariant under the group (\ref{affine}).  It can be shown similarly that each quadratic form $(\bY_j - \bmuhat_2)'\bSigmahat^{-1}_{22}(\bY_j-\bmuhat_2)$, $j=n+1,\ldots,N$ is invariant under (\ref{affine}), hence so is $b_{2,p,q}^{(2)}$.  Therefore, in the multivariate normal case, the invariance of $b_{2,p,q}$ under (\ref{affine}) is a consequence of equivariance.

In non-normal cases, however, the detailed computations in the proof of Theorem \ref{invariant} are necessary to prove that $b_{2,p,q}$ is invariant under (\ref{affine}).  Moreover, now that it has been established that $b_{2,p,q}$ is invariant under (\ref{affine}), we then choose $\bLambda_{11} = \bSigma_{11 \cdot 2}^{-1/2}$, $\bLambda_{22} = \bSigma_{22}^{-1/2}$, $\bLambda_{12} = -\bSigma_{12}\bSigma_{22}^{-1}$, and $\bnu = - \bLambda\bC\bmu$ to reduce the data to being mutually independent and monotone incomplete with mean $\bzero$ and covariance matrix $\bLambda\bC\bSigma\bC'\bLambda' = \bI_{p+q}$.  Therefore, in deriving the distribution of $b_{2,p,q}$, we assume without loss of generality that $\bmu = \bzero$ and $\bSigma = \bI_{p+q}$.  
}\end{remark}

\section{The null and non-null asymptotic distributions of \texorpdfstring{$\boldsymbol{b_{2,p,q}}$}{b2pq}}
\label{asymptoticdistns}

For the $(p+q)$-dimensional random vector $\bZ = \begin{pmatrix} \bX \\ \bY \end{pmatrix}$, define 
\begin{equation}
\label{Xidef}
\bXi = E\big((\bZ\bZ')^2\big) = E\big(\bZ(\bZ'\bZ)\bZ'\big) = E(\|\bZ\|^2\,\bZ\bZ'),
\end{equation}
and write $\bXi$ in partitioned form, 
\begin{equation}
\label{Xi}
\bXi \equiv \begin{pmatrix}\bXi_{11}&\bXi_{12}\\\bXi_{21}&\bXi_{22}\end{pmatrix} = \begin{pmatrix}E(\|\bZ\|^2\,\bX\bX')&E(\|\bZ\|^2\,\bX\bY')\\E(\|\bZ\|^2\,\bY\bX')&E(\|\bZ\|^2\,\bY\bY')\end{pmatrix}.
\end{equation}
Define 
\begin{equation}
\label{Xistar}
\bXi^* = E\big((\bY\bY')^2\big), \qquad \bTheta^*=E(\|\bY\|^2\,\bY);
\end{equation}
and set 
\begin{equation}
\label{Xitilde}
\wtilde\bXi \equiv 
\begin{pmatrix}
\wtilde\bXi_{11}&\wtilde\bXi_{12}\\
\wtilde\bXi_{21}&\wtilde\bXi_{22}
\end{pmatrix} =
\begin{pmatrix}
c_1\bXi_{11}&c_1\bXi_{12}\\
c_1\bXi_{21}&c_1\tau\bXi_{22}+c_2\bartau\bXi^*
\end{pmatrix},
\end{equation}
\begin{equation}
\label{Theta}
\bTheta \equiv \begin{pmatrix}\bTheta_1 \\ \bTheta_2\end{pmatrix} = 
\begin{pmatrix}E(\|\bZ\|^2\bX) \\ E(\|\bZ\|^2\bY) \end{pmatrix},
\end{equation}
and 
\begin{equation}
\label{Thetatilde}
\wtilde\bTheta \equiv \begin{pmatrix}\wtilde\bTheta_1 \\ \wtilde\bTheta_2 \end{pmatrix} = \begin{pmatrix}c_1\bTheta_1 \\ c_1\tau\bTheta_2+c_2\bartau\bTheta^*\end{pmatrix}.
\end{equation}
Our main result, given in the following theorem, provides the non-null distribution of $b_{2,p,q}$ for a large class of alternatives.  

\begin{theorem}
\label{basympthm}
Suppose that the monotone incomplete random sample (\ref{monotonesample}) is drawn from a population modeled by a random vector $\bZ=\begin{pmatrix}\bX\\\bY\end{pmatrix}$ such that $E(\bZ) = \bzero$, $\Cov(\bZ) = \bI_{p+q}$, and $E\|\bZ\|^8 < \infty$.  For $n, N \rightarrow \infty$ with $n/N \to \delta \in (0,1)$, we have 
\begin{equation}
\label{nonnulldistn}
N^{1/2}\(b_{2,p,q} - \nu\)/\sigma \toD N(0,1),
\end{equation}
where 
\begin{equation}
\label{nu}
\nu = c_1\tau E(\|\bZ\|^4) + c_2\bartau E(\|\bY\|^4)
\end{equation}
and 
\begin{equation}
\label{sigmasq}
\begin{aligned}
\sigma^2 \ = \ & \tau \Big\{c_1^2 \, \Var\|\bZ\|^4 + 4 \, \Var(\bZ'\widetilde\bXi\bZ) + 16 \, \widetilde\bTheta'\widetilde\bTheta \\
& \qquad 
- 4 \, c_1 \, \Cov(\|\bZ\|^4, \bZ'\widetilde\bXi\bZ) - 8 c_1 \, E\|\bZ\|^4 \bZ'\widetilde\bTheta + 16 \, E[\bZ'\widetilde\bXi\bZ\bZ']\widetilde\bTheta\Big\} \\
& + \bartau\Big\{c_2^2 \, \Var\|\bY\|^4 + 4 \, \Var(\bY'\widetilde\bXi_{22}\bY) + 16 \, \widetilde\bTheta_2'\widetilde\bTheta_2 \\
& \qquad
-4 \, c_2 \Cov(\|\bY\|^4, \bY'\widetilde\bXi_{22}\bY)
-8 c_2 \, E\|\bY\|^4\bY'\widetilde\bTheta_2 + 16 \, E[\bY'\widetilde\bXi_{22}\bY\bY']\widetilde\bTheta_2\Big\}.
\end{aligned}
\end{equation}
\end{theorem}

The proof of this result is provided in Appendix \ref{basymptotics}.  For the case in which $\bZ \sim N_{p+q}(\bzero,\bI_{p+q})$, the limiting distribution (\ref{nonnulldistn}) reduces to the following result on the null distribution of $b_{2,p,q}$.  

\begin{corollary}
\label{basympcor}
Suppose that the monotone incomplete  sample (\ref{monotonesample}) is drawn from $N_{p+q}(\bzero,\bI_{p+q})$.  For $n, N \rightarrow \infty$ with $n/N \to \delta \in (0,1)$, we have 
\begin{equation}
\label{nulldistn}
N^{1/2}\(b_{2,p,q} - \nu\)/\sigma \toD N(0,1),
\end{equation}
where 
\begin{equation}
\label{nullnu}
\nu = c_1\tau (p+q)(p+q+2) + c_2\bartau q(q+2)
\end{equation}
and 
\begin{equation}
\label{nullsigmasq}
\begin{aligned}
\sigma^2 = \ & 8 \tau\Big\{c^2_1(p+q)(p+q+2)(p+q+3)+ c^2_1p(p+q+2)^2 \\
& \qquad + q\big(c_1\tau(p+q+2)+c_2\bar\tau(q+2)\big)^2 \\
& \qquad - 2c_1(p+q+2)\Big(pc_1(p+q+2) + q\big(c_1\tau(p+q+2) + c_2\bar\tau(q+2)\big)\Big)\Big\} \\
& + 8 \bar\tau\Big\{c_2^2q(q+2)(q+3)+ q\big(c_1\tau(p+q+2)+c_2\bar\tau(q+2)\big)^2 \\
& \qquad\quad - 2c_2(q+2)q\big(c_1\tau (p+q+2)+ c_2\bar\tau(q+2)\big)\Big\}.
\end{aligned}
\end{equation}
\end{corollary}

\begin{remark}
\label{specialcases}
{\rm 
For $(c_1,c_2) = (\tau,\bartau)$, (\ref{nullnu}) and (\ref{nullsigmasq}) reduce to 
\begin{align}
\label{taumean}
\nu & = \tau^2 (p+q)(p+q+2) + \bar\tau^2 q(q+2) \\
\intertext{and} 
\label{tauvariance}
\sigma^2 & = 8\{\tau^3 (p+q)(p+q+2) + \bar\tau^3 q(q+2) + \tau\bar\tau(\tau(p+q+2)-\bar\tau(q+2))^2\},
\end{align}
respectively.  For $(c_1,c_2) = (1,1)$, (\ref{nullnu}) and (\ref{nullsigmasq}) reduce, respectively, to 
\begin{align}
\label{mardiamean}
\nu & = \tau (p + q)(p + q + 2) + \bar\tau q(q + 2) \\
\intertext{and}
\label{mardiavariance}
\sigma^2 & = 8\{\tau (p + q)(p + q + 2) + \bar\tau q(q + 2)+\tau\bar\tau p^2q\}.
\end{align}
Corollary \ref{basympcor} also generalizes the result of Mardia (1970) for complete data; indeed, when we set $\tau =1$, i.e., $n = N$ in (\ref{mardiamean}) and (\ref{mardiavariance}), we find that $\nu$ reduces to $(p+q)(p+q+2)$ and $\sigma^2$ reduces to $8(p+q)(p+q+2)$; hence (\ref{nulldistn}) reduces to the result of Mardia.  
}\end{remark}

\subsection{Application to the Pennsylvania cholesterol data set}

Assume that the Pennsylvania cholesterol data (Ryan, et al. 2005, p. 267) consist of mutually independent vectors and that missing observations are MCAR.  For that data set, we have $p = 1$, $q = 2$, $N = 28$, and $n = 19$, and we choose $c_1 = \tau$ and $c_2 = \bartau$.   The asymptotic mean and variance of $b_{2,1,2}$ are obtained from (\ref{nullnu}) and (\ref{nullsigmasq}) to be $\nu = 7.7334$ and $\sigma^2 = 181.1658$, respectively.  By Corollary \ref{basympcor}, the asymptotic null distribution of the statistic (\ref{nulldistn}) is 
$$
\sqrt{28}(b_{2,1,2} - 7.7334)/13.4598 \approx N(0,1).
$$
We calculate using (\ref{b-formula}) that the observed value of $b_{2,1,2}$ is $5.8623$.  Therefore, the observed value of the statistic (\ref{nulldistn}) is 
$$
\sqrt{28}(5.8623 - 7.7334)/13.4598 = -0.7356.
$$
The approximate $P$-value of the test is $2\Phi(-0.7356) = 0.2310$, where $\Phi(\cdot)$ denotes the cumulative distribution function of the standard normal distribution.  Therefore, we fail to reject the null hypothesis of multivariate normality at the 5\% level of significance.

We note that the same conclusion is obtained by applying the classical Mardia statistic to the subset of the Pennsylvania cholesterol data set consisting of the $n = 19$ complete observations only.  For this subset, the observed value of Mardia's statistic is $7.8176$, the corresponding observed value of the normal approximation to Mardia's statistic is $-0.1207$, hence the resulting approximate $P$-value for Mardia test is $2\Phi(-0.1207) = 0.9038$.  However, this $P$-value is so large that the test based on the complete data appears unable to assess the strength of the evidence against the null hypothesis of normality, and this reflects the loss of information inherent in discarding the incomplete observations; cf., Little and Rubin (2002), p. 41.  By contrast, the $P$-value based on the full data set appears to provide some measure of the strength of the evidence against the null hypothesis, even though the strength of that evidence is assessed to be too weak to reject that hypothesis.  
  
We remark also that, in the case of the cholesterol data, the smaller sample size is less likely to yield an accurate normal approximation to Mardia's statistic, so the substantially larger $P$-value of the Mardia test for the complete data perhaps should be applied cautiously.

\vskip 0.4truein

\noindent
{\bf Acknowledgments.}  We are grateful to the referees and the associate editor for comments which helped us to improve the manuscript.  The research of Richards was also supported by a 2013--2014 sabbatical leave-of-absence and a Romberg Guest Professorship at the Heidelberg University Graduate School for Mathematical and Computational Methods in the Sciences, funded by the German Universities Excellence Initiative grant GSC 220/2.

\bigskip
\bigskip
\bigskip

{
\noindent
{\bf\Large References}
\medskip
\parskip=3pt
\parindent=0pt

Anderson, T.~W.~(2003).  {\sl An Introduction to Multivariate Statistical Analysis} (third edition).  Wiley, New York.

Chang, W.-Y., and Richards, D.~St.~P. (2009).  Finite-sample inference with monotone incomplete multivariate normal data, I. {\it J. Multivariate Anal.}, {\bf 100}, 1883--1899.

Chang, W.-Y., and Richards, D.~St.~P. (2010).  Finite-sample inference with monotone incomplete multivariate normal data, II. {\it J. Multivariate Anal.}, {\bf 101}, 603--620.

Cohn, N., Davidov, O., and Haitovsky, Y. (2008).  Double sampling designs in multivariate linear models with missing data.  {\it Comm. Statist. Simulation Comput.}, {\bf 37}, 1156--1166.

Davidov, O., and Peddada, S.~D. (2013).  The linear stochastic order and directed inference for multivariate ordered distributions.  {\it Ann. Statist.}, {\bf 41}, 1--40.

Eaton, M. L., and Kariya, T. (1983).  Multivariate tests with incomplete data.  {\it Ann. Statist.}, {\bf 11}, 654--665.

Garren, S. T., and Peddada, S. D. (2000).  Asymptotic normality in multivariate nonlinear regression and multivariate generalized linear regression models under repeated measurements with missing data. {\it Statist. Probab. Lett.}, {\bf 48}, 293--302. 

Hao, J. and Krishnamoorthy, K. (2001).  Inferences on a normal covariance matrix and generalized variance with monotone missing data.  {\it J. Multivariate Anal.}, {\bf 78}, 62--82.

Henze, N.~(1994).  On Mardia's kurtosis test for multivariate normality.  {\it Commun. Statist. -- Theory \& Methods}, {\bf 23}, 1031--1045.

Henze, N. (2002).  Invariant tests for multivariate normality: A critical review.  {\it Statist. Papers}, {\bf 43}, 467--506.

Krishnamoorthy, K., and Yu, J. (2012).  Multivariate Behrens-Fisher problem with missing data. {\it J. Multivariate Anal.}, {\bf 105}, 141--150.

Little, R. J. A. (1976).  Inference about means from incomplete multivariate data, {\it Biometrika} 63 (1976) 593--604.

Little,  R. J. A., and Rubin, D. B. (2002).  {\sl Statistical Analysis with Missing Data} (second edition).  Wiley, Hoboken, NJ.

Mardia, K.~V. (1970).  Measures of multivariate skewness and kurtosis with applications.  {\it Biometrika}, {\bf 57}, 519--530. 

Mardia, K.~V. (1974).  Applications of some measures of multivariate skewness and kurtosis in testing normality and robustness studies.  {\it Sankhy\=a} B, {\bf 36}, 115--128. 

Peddada, S.~D., Harris, S., Davidov, O. (2010).  Analysis of correlated gene expression data on ordered categories.  {\it J. Ind. Soc. Agric. Statist.}, {\bf 64}, 45--60.

Richards, D.~St.~P., and Yamada, T. (2010).  The Stein phenomenon for monotone incomplete multivariate normal data.  {\it J. Multivariate Anal.}, {\bf 101}, 657--678.

Romer, M.~M. (2009). {\sl The Statistical Analysis of Monotone Incomplete Multivariate Normal Data}.  Doctoral Dissertation, Penn State University.

Romer, M. M., and Richards, D.~St.~P. (2010).  Maximum likelihood estimation of the mean of a multivariate normal population with monotone incomplete data.  {\it Statist. \& Probab. Lett.}, {\bf 80}, 1284--1288.

Romer, M.~M., and Richards, D.~St.~P. (2013).  Finite-sample inference with monotone incomplete multivariate normal data, III: Hotelling's $T^2$-statistic.  {\it Statist. Modelling}, {\bf 13}, 431--457.

Ryan, B., Joiner, B., and Cryer, J. (2005).  {\sl Minitab Handbook} (fifth edition).  Duxbury Press, Boston.

Yamada, T. (2013).  Asymptotic properties of canonical correlation analysis for one group with additional observations.  {\it J. Multivariate Anal.}, {\bf 114}, 389--401.
}

\appendix

\section{Asymptotic expansions of 
\texorpdfstring{$\bSigmahat$}{Sigmahat} and 
\texorpdfstring{$\bSigmahat^{-1}$}{Sigmahatinverse}}
\label{Sigma-asymptotics}
\setcounter{equation}{0}

To derive the expansions, we utilize results of Henze (1994).  For a sequence of random vectors $\bomega_n$ and positive scalars $a_n$, we use the notation $\bomega_n = O_p(a_n)$ to denote that the sequence $\bomega_n/a_n$ is {\it tight}, and we write $\bomega_n = o(a_n)$ to mean that $\bomega_n/a_n$ {\it converges to zero in probability}.  

Let $\bS = n^{-1} \bA$, where $\bA$ is given in (\ref{Amatrix}).  Since $\bA$ is the matrix of squares and cross-products of the $n$ complete observations in the data set (\ref{monotonesample}) then, following Henze (1994, p. 1035), we deduce that 
$\bS = \bI_{p+q} + n^{-1/2}\bB + O_p(n^{-1})$, where 
$$
\bB = n^{1/2} \Big(\frac1n\sum_{j=1}^n \bZ_j\bZ_j' - \bI_{p+q}\Big) = O_p(1).
$$
Write $\bB = \begin{pmatrix}\bB_{11}&\bB_{12}\\\bB_{21}&\bB_{22}\end{pmatrix}$, where 
\begin{align*}
\bB_{11} &= n^{1/2} \Big(\dfrac1n\sum\limits_{j=1}^n\bX_j\bX_j'-\bI_{p}\Big), \\
\bB_{12} &= \bB_{21}'= n ^{-1/2} \sum\limits_{j=1}^n\bX_j\bY_j', \\
\bB_{22} &= n^{1/2} \Big(\frac1n\sum_{j=1}^n\bY_j\bY_j'-\bI_{q}\Big),
\end{align*}
and define $\bB^* = \tau^{-1/2}\bB = 
\begin{pmatrix} \bB_{11}^* & \bB_{12}^* \\
\bB_{21}^* & \bB_{22}^* \end{pmatrix}$.  
Then, 
$$
\bA = \begin{pmatrix} \bA_{11} & \bA_{12} \\
\bA_{21} & \bA_{22,n} \end{pmatrix} = n\bI_{p+q}+n^{1/2}\bB+O_p(1),
$$
so that 
\begin{equation}
\label{Aexpansion}
\begin{aligned}
\bA_{11}&=n\bI_{p}+n^{1/2}\bB_{11} +O_p(1), \\
\bA_{12}&=\bA_{21}'= n^{1/2}\bB_{12} +O_p(1), \\
\bA_{22,n}&=n\bI_{q}+n^{1/2}\bB_{22} +O_p(1).
\end{aligned}
\end{equation}
Next, we follow Henze (1994, p. 1035, equation (2.2)) to obtain 
$$
\bA_{22,n}^{-1} = n^{-1}\bI_{q}-n^{-3/2}\bB_{22} +O_p(n^{-2}).
$$
By the same argument, we obtain 
\begin{align*}
\bA_{22,N}&=N\bI_{q}+N^{1/2}\wtilde\bB_{22} +O_p(1), \\
\intertext{and, by inversion, } 
\bA_{22,N}^{-1}&=N^{-1}\bI_{q}-N^{-3/2}\wtilde\bB_{22} +O_p(N^{-2}), 
\end{align*} 
where 
$$
\wtilde\bB_{22} = N^{1/2} \Big(\frac1N \sum_{j=1}^N \bY_j\bY_j'-\bI_q\Big) 
= O_p(1).
$$

The {\it partial Iwasawa coordinates} of $\bSigmahat$ (Chang and Richards (2009, equation (4.5)) are 
\begin{equation}
\label{Deltahatij}
\bDeltahat_{11} = \frac1n \bA_{11\cdot2,n}, \ \ \bDeltahat_{12}=\bDeltahat_{21}'=\bA_{12}\bA_{22,n}^{-1}, \ \ \bDeltahat_{22}=\frac1N\bA_{22,N}.
\end{equation}
These matrices are such that 
\begin{align}
\label{iwasawa}
\bSigmahat&=
\begin{pmatrix} \bI_p & \bDeltahat_{12} \\ \bzero & \bI_q \end{pmatrix}
\begin{pmatrix} \bDeltahat_{11} & \bzero \\ \bzero & \bDeltahat_{22} \end{pmatrix}
\begin{pmatrix} \bI_p & \bzero \\ \bDeltahat_{21} & \bI_q \end{pmatrix},
\end{align}
hence 
\begin{align}
\label{Sigmahatinv}
\bSigmahat^{-1} & = 
\begin{pmatrix}
\bDeltahat_{11}^{-1}&-\bDeltahat_{11}^{-1}\bDeltahat_{12}\\
-\bDeltahat_{21}\bDeltahat_{11}^{-1}&
\bDeltahat_{22}^{-1}+\bDeltahat_{21}\bDeltahat_{11}^{-1}\bDeltahat_{12}
\end{pmatrix}.
\end{align}
Noting that $n=\tau N$ and $\bB^*=\tau^{-1/2}\bB$, we apply (\ref{Aexpansion}) to deduce that $\bDeltahat_{11}$, $\bDeltahat_{12}$, and $\bDeltahat_{22}$ have the following asymptotic expansions:
\begin{align*}
\bDeltahat_{11} & = \frac1n\{\bA_{11}-\bA_{12}\bA_{22,n}^{-1}\bA_{21}\}\\
& = \frac1n \{n\bI_{p}+n^{1/2}\bB_{11}-(n^{1/2}\bB_{12})(n^{-1}\bI_{q}-n^{-3/2}\bB_{22})(n^{1/2}\bB_{21})\} +O_p(n^{-1})\\
& = \bI_{p}+n^{-1/2}\bB_{11}+O_p(n^{-1}) \\
& = \bI_{p}+N^{-1/2}\bB_{11}^*+O_p(N^{-1});
\end{align*}
\begin{align*}
\bDeltahat_{12} & = \bA_{12}\bA_{22,n}^{-1}=n^{1/2}\bB_{12}(n^{-1}\bI_{q}-n^{-3/2}\bB_{22})+O_p(n^{-1})\\
& = n^{-1/2}\bB_{12}-n^{-1}\bB_{12}\bB_{22}+O_p(n^{-1}) \\
& = N^{-1/2}\bB_{12}^*-N^{-1}\bB_{12}^*\bB_{22}^*+O_p(N^{-1});
\intertext{and}
\bDeltahat_{22} & = \frac1N\bA_{22,N} 
= \bI_{q}+N^{-1/2}\wtilde\bB_{22} +O_p(N^{-1}).
\end{align*} 
By the same inversion argument used earlier, we also have 
\begin{align*}
\bDeltahat_{11}^{-1} = \bI_{p}-N^{-1/2}\bB_{11}^*+O_p(N^{-1}), \qquad
\bDeltahat_{22}^{-1} = \bI_{q}-N^{-1/2}\wtilde\bB_{22} +O_p(N^{-1}).
\end{align*} 
Collecting together these results, we obtain 
\begin{align*}
\bSigmahat_{11} &= \bDeltahat_{11} + \bDeltahat_{12}\bDeltahat_{22}^{-1}\bDeltahat_{21} = \bI_{p}+N^{-1/2}\bB_{11}^* + O_p(N^{-1}), \\
\bSigmahat_{12} &= \bDeltahat_{12}\bDeltahat_{22} = N^{-1/2}\bB_{12}^* + O_p(N^{-1}), \\
\bSigmahat_{22} &= \bDeltahat_{22} = \bI_{q}+N^{-1/2}\wtilde\bB_{22} + O_p(N^{-1}).
\end{align*}
Therefore, 
\begin{align}
\label{Sigmahatasymp}
\bSigmahat = \bI_{p+q} + N^{-1/2}\wtilde\bB + O_p(N^{-1}),
\end{align}
where $\wtilde\bB = \begin{pmatrix} 
\bB_{11}^* & \bB_{12}^* \\ \bB_{21}^*&\wtilde\bB_{22}
\end{pmatrix}$. 
By inverting (\ref{Sigmahatasymp}), we obtain the asymptotic expansion, 
\begin{align}
\bSigmahat^{-1} = \bI_{p+q}-N^{-1/2}\wtilde\bB + O_p(N^{-1}).
\end{align}

\section{Asymptotic expansions of \texorpdfstring{$\boldsymbol{b_{2,p,q}^{(1)}}$}{b1} and \texorpdfstring{$\boldsymbol{b_{2,p,q}^{(2)}}$}{b2}}
\label{b1b2asymptotics}
\setcounter{equation}{0}

The statistic $b_{2,p,q}^{(2)}$ is reminiscent of a statistic studied by Henze (1994) so we shall proceed similarly to Henze in obtaining an asymptotic expansion for $b_{2,p,q}^{(2)}$.  Nevertheless, we note that the statistic $b_{2,p,q}^{(2)}$ is not identical with Henze's statistic; indeed, although $b_{2,p,q}^{(2)}$ is given by a sum over the indices $j = n+1,\ldots,N$, we note that $\bar\bY$ and $\bSigmahat_{22} = N^{-1}\bA_{22,N}$ are the sample mean and sample covariance matrix, respectively, of $\bY_1,\ldots,\bY_N$, so that $b_{2,p,q}^{(2)}$ incorporates information from the portion of the data that is complete.

Thus, by direct expansion, we have
\begin{align}
\label{yjminusybar}
[(\bY_j-\bar\bY)'\bA_{22,N}^{-1}(\bY_j-\bar\bY)]^2 = \ & (\bY_j'\bA_{22,N}^{-1}\bY_j - 2\bar\bY'\bA_{22,N}^{-1}\bY_j + \bar\bY'\bA_{22,N}^{-1}\bar\bY )^2 \nonumber \\
= \ & (\bY_j'\bA_{22,N}^{-1}\bY_j)^2 + 4(\bar\bY'\bA_{22,N}^{-1}\bY_j)^2 \nonumber \\
& + (\bar\bY'\bA_{22,N}^{-1}\bar\bY)^2 - 4\bY_j'\bA_{22,N}^{-1}\bY_j\bar\bY'\bA_{22,N}^{-1}\bY_j \\
& + 2\bY_j'\bA_{22,N}^{-1}\bY_j\bar\bY'\bA_{22,N}^{-1}\bar\bY 
- 4\bar\bY'\bA_{22,N}^{-1}\bY_j\bar\bY'\bA_{22,N}^{-1}\bar\bY. \qquad \nonumber 
\end{align}
With $\bXi^*$ and $\bTheta^*$ defined as in (\ref{Xistar}) we obtain, by applying the methods of Henze (1994), the approximations  
$$
\frac1N \sum_{j=n+1}^N (\bY_j\bY_j')^2 = \bartau\bXi^* + o(1)
$$
and 
$$
\frac1N \sum_{j=n+1}^N \|\bY_j\|^2\,\bY_j = \bartau\bTheta^* + o(1).
$$
Therefore, 
$$
\frac1N \sum_{j=n+1}^N (\bY_j'\bA_{22,N}^{-1}\bY_j)^2 = 
\frac1N \sum_{j=n+1}^N(\bY_j'\bY_j)^2 
- 2N^{-1/2} \bartau \, \tr(\wtilde\bB_{22}\bXi^*) + o(N^{-1/2}),
$$
and 
$$
\frac1N \sum_{j=n+1}^N \bY_j'\bA_{22,N}^{-1}\bY_j\bar\bY'\bA_{22,N}^{-1}\bY_j = \bartau{\bTheta^*}'\bar\bY+o(N^{-1/2}).
$$
Also, 
$$
N^{-1}\sum_{j=n+1}^N(\bar\bY '\bA_{22,N}^{-1}\bY_j )^2 = O_p(N^{-1}),
$$
$$
N^{-1}\sum_{j=n+1}^N(\bar\bY '\bA_{22,N}^{-1}\bar\bY )^2 = O_p(N^{-2}),
$$
$$
N^{-1}\sum_{j=n+1}^N\bY_j'\bA_{22,N}^{-1}\bY_j\bar\bY'\bA_{22,N}^{-1}\bar\bY = O_p(N^{-1}),
$$
and 
$$
N^{-1}\sum_{j=n+1}^N\bar\bY'\bA_{22,N}^{-1}\bY_j\bar\bY'\bA_{22,N}^{-1}\bar\bY = O_p(N^{-1}).
$$
Summing (\ref{yjminusybar}) over $j=n+1,\ldots,N$, we obtain 
\begin{equation}
\label{asymptT2}
\frac1N b_{2,p,q}^{(2)} = \frac1N \sum_{j=n+1}^N (\bY_j'\bY_j)^2 - 2N^{-1/2}\bartau \, \tr(\wtilde\bB_{22}\bXi^*) - 4\bartau{\bTheta^*}'\bar\bY + o(N^{-1/2}).
\end{equation}

Consider next the statistic $b_{2,p,q}^{(1)}$.  We apply (\ref{zjminusmuhat}) to write 
\begin{align*}
[(\bZ_j-\bmuhat)'\bSigmahat^{-1}(\bZ_j-\bmuhat)]^2 & =
[(\bZ_j-\bar\bZ_1 + \bartau\wtilde\bY)'\bSigmahat^{-1}(\bZ_j-\bar\bZ_1 + \bartau\wtilde\bY)]^2 \\
& = [(\bZ_j-\bar\bZ_1)'\bSigmahat^{-1}(\bZ_j-\bar\bZ_1) + 2\bartau(\bZ_j-\bar\bZ_1)'\bSigmahat^{-1}\wtilde\bY + \bartau^2\wtilde\bY'\bSigmahat^{-1}\wtilde\bY]^2,
\end{align*}
expand this expression directly, and sum over $j=1,\ldots,n$.  By (\ref{z1bar}), $\sum_{j=1}^n (\bZ_j-\bar\bZ_1) = \bzero$, so we obtain $b_{2,p,q}^{(1)}$ as a sum of five terms, $b_{2,p,q}^{(1)} = \sum_{j=1}^5 T_j$, where 
\begin{eqnarray*}
T_1 & = & \sum_{j=1}^n [(\bZ_j-\bar\bZ_1)'\bSigmahat^{-1}(\bZ_j-\bar\bZ_1)]^2, \\
T_2 & = & 4\bartau\sum_{j=1}^n (\bZ_j-\bar\bZ_1)'\bSigmahat^{-1}(\bZ_j-\bar\bZ_1)(\bZ_j-\bar\bZ_1)'\bSigmahat^{-1}\wtilde\bY, \\
T_3 & = & 4\bartau^2\sum_{j=1}^n (\bZ_j-\bar\bZ_1)'\bSigmahat^{-1}\wtilde\bY{\wtilde\bY}'\bSigmahat^{-1}(\bZ_j-\bar\bZ_1), \\
T_4 & = & \bartau^2 \sum_{j=1}^n (\bZ_j-\bar\bZ_1)'\bSigmahat^{-1}(\bZ_j-\bar\bZ_1){\wtilde\bY}'\bSigmahat^{-1}\wtilde\bY, \\
T_5 & = & n\bartau^2({\wtilde\bY}'\bSigmahat^{-1}\wtilde\bY)^2.
\end{eqnarray*}
From (\ref{Amatrix}), 
\begin{align*}
\sum_{j=1}^n(\bZ_j-\bar\bZ_1)(\bZ_j-\bar\bZ_1)' = \bA \equiv 
\begin{pmatrix}\bA_{11} &\bA_{12}\\ \bA_{21}&\bA_{22,n}\end{pmatrix}
\end{align*}
and by (\ref{Sigmahatinv}), 
$$
\bSigmahat^{-1} 
= \begin{pmatrix}
n\bA_{11\cdot2}^{-1} &-n\bA_{11\cdot2}^{-1}\bA_{12}\bA_{22,n}^{-1}\\
-n\bA_{22,n}^{-1}\bA_{21}\bA_{11\cdot2}^{-1}&
N\bA_{22,N}^{-1}+n\bA_{22,n}^{-1}\bA_{21}\bA_{11\cdot2}^{-1}\bA_{12}\bA_{22,n}^{-1}
\end{pmatrix}. 
$$
By direct multiplication, we obtain 
\begin{equation}
\label{SigmahatinvA}
\bSigmahat^{-1}\bA = \begin{pmatrix}n\bI_p &\bzero\\ n\bA_{22,n}^{-1}\bA_{21}+N\bA_{22,N}^{-1}\bA_{21}&N\bA_{22,N}^{-1}\bA_{22,n}\end{pmatrix}
\end{equation}
and 
\begin{equation}
\label{ASigmahatinvA}
\bA\bSigmahat^{-1}\bA = \begin{pmatrix}n\bA_{11\cdot 2,n}+N \bA_{12}\bA_{22,N}^{-1}\bA_{21}&N\bA_{12}\bA_{22,N}^{-1}\bA_{22,n}\\ N\bA_{22,n}\bA_{22,N}^{-1}\bA_{21}&N\bA_{22,n}\bA_{22,N}^{-1}\bA_{22,n}\end{pmatrix}, \\
\end{equation}
and then it follows from (\ref{ytilde}) that 
\begin{eqnarray*}
\wtilde\bY\bSigmahat^{-1}\wtilde\bY & = & \begin{pmatrix} \bzero \\ \bar\bY_1-\bar\bY_2 \end{pmatrix}' \begin{pmatrix} \bzero &\bzero\\ \bzero&\bA_{22,n}^{-1} \end{pmatrix}\bA\bSigmahat^{-1}\bA \begin{pmatrix} \bzero &\bzero \\ \bzero&\bA_{22,n}^{-1}\end{pmatrix} \begin{pmatrix}\bzero\\ \bar\bY_1-\bar \bY_2 
\end{pmatrix} \\
& = & N(\bar\bY_1-\bar\bY_2)'\bA_{22,N}^{-1}(\bar\bY_1-\bar\bY_2).
\end{eqnarray*}
Also, $N^{1/2}(\bar\bY_1-\bar\bY_2)(\bar\bY_1-\bar\bY_2)'=O_p(N^{-1/2})$.

The statistic $T_1$ is the same as in the complete case, so we apply the expansion derived by Henze (1994).  Letting $\bXi = E\((\bZ\bZ')^2\)$ and $\bTheta = E(\|\bZ\|^2\,\bZ)$, we obtain 
\begin{align*}
\frac1N T_1 = &
\frac1N \sum_{j=1}^n [(\bZ_j-\bar\bZ_1)'\bSigmahat^{-1}(\bZ_j-\bar\bZ_1)]^2 \\ 
= & \frac1N\sum_{j=1}^n(\bZ_j '\bZ_j )^2-2N^{-1/2}\tau \, \tr(\wtilde\bB\bXi )-4\tau\bTheta'\bar\bZ_1 +o(N^{-1/2}). 
\end{align*}

We now consider $T_2$.  Letting $\wtilde\bY^*=\begin{pmatrix}\bzero\\ \bar \bY_1-\bar \bY_2 \end{pmatrix}$, we obtain 
$$
\wtilde\bY = \begin{pmatrix}\bzero & \bA_{12}\bA_{22,n}^{-1} \\ \bzero&\bI_q\end{pmatrix}\wtilde\bY^* = \begin{pmatrix}\bzero &\bzero\\ \bzero&\bI_q\end{pmatrix}\wtilde\bY^* 
+N^{-1/2}\begin{pmatrix}\bzero &\bB_{12}\\ \bzero&\bI_q\end{pmatrix}\wtilde\bY^* + O_p(N^{-1}),
$$
and 
\begin{multline*}
(\bZ_j-\bar\bZ_1)'\bSigmahat^{-1}(\bZ_j-\bar\bZ_1)(\bZ_j-\bar\bZ_1)'\bSigmahat^{-1}\wtilde\bY\\
=\bZ_j'\bSigmahat^{-1}\bZ_j\bZ_j'\bSigmahat^{-1}\wtilde\bY
-\bZ_j'\bSigmahat^{-1}\bZ_j\bar \bZ'\bSigmahat^{-1}\wtilde\bY
-2\bZ_j'\bSigmahat^{-1}\bar \bZ\bZ_j'\bSigmahat^{-1}\wtilde\bY\\
+2\bZ_j'\bSigmahat^{-1}\bar \bZ\bar \bZ'\bSigmahat^{-1}\wtilde\bY
+\bar \bZ'\bSigmahat^{-1}\bar \bZ\bZ_j'\bSigmahat^{-1}\wtilde\bY
-\bar \bZ'\bSigmahat^{-1}\bar \bZ\bar \bZ'\bSigmahat^{-1}\wtilde\bY.
\end{multline*}
Then,
\begin{align*}
\frac1N\sum_{j=1}^n&\bZ_j'\bSigmahat^{-1}\bZ_j\bZ_j'\bSigmahat^{-1}\wtilde\bY = \tau\bTheta'\wtilde\bY^* +o(N^{-1/2}),\\
\frac1N\sum_{j=1}^n&\bZ_j'\bSigmahat^{-1}\bZ_j\bar \bZ'\bSigmahat^{-1}\wtilde\bY = O_p(N^{-1}),\\
\frac1N\sum_{j=1}^n&\bZ_j'\bSigmahat^{-1}\bar \bZ\bZ_j'\bSigmahat^{-1}\wtilde\bY = O_p(N^{-1}),
\end{align*}
and
\begin{multline*}
\frac1N\sum_{j=1}^n 
\{2\bZ_j'\bSigmahat^{-1}\bar \bZ\bar \bZ'\bSigmahat^{-1}\wtilde\bY
+\bar \bZ'\bSigmahat^{-1}\bar \bZ\bZ_j'\bSigmahat^{-1}\wtilde\bY
-\bar \bZ'\bSigmahat^{-1}\bar \bZ\bar \bZ'\bSigmahat^{-1}\wtilde\bY\} \\ 
= 2\tau\bar \bZ'\bSigmahat^{-1}\bar \bZ\bar \bZ'\bSigmahat^{-1}\wtilde\bY = O_p(N^{-2}).
\end{multline*}
Therefore, 
\begin{align*}
\frac1{4\bartau N} T_2 = & 
\frac1N\sum_{j=1}^n(\bZ_j-\bar\bZ_1)'\bSigmahat^{-1}(\bZ_j-\bar\bZ_1)(\bZ_j-\bar\bZ_1)'\bSigmahat^{-1}\wtilde\bY \\
= & \tau\bTheta'\wtilde\bY^* +o(N^{-1/2}).
\end{align*}

Next, 
\begin{align*}
\frac1{4\bartau^2}T_3 & = \sum_{j=1}^n(\bZ_j-\bar\bZ_1)'\bSigmahat^{-1}\wtilde\bY{\wtilde\bY}'\bSigmahat^{-1}(\bZ_j-\bar\bZ_1)\\&={\wtilde\bY}'\bSigmahat^{-1}\Big[\sum_{j=1}^n(\bZ_j-\bar\bZ_1)(\bZ_j-\bar\bZ_1)'\Big]\bSigmahat^{-1}\wtilde\bY \\
& = {\wtilde\bY}'\bSigmahat^{-1}\bA\bSigmahat^{-1}\wtilde\bY\\
& =
\begin{pmatrix}
\bzero\\\bar  \bY_1-\bar \bY_2 
\end{pmatrix}'
\begin{pmatrix}
\bzero &\bzero\\ \bzero&\bA_{22,n}^{-1}
\end{pmatrix}
\bA\bSigmahat^{-1}\bA\bSigmahat^{-1}\bA
\begin{pmatrix}
\bzero &\bzero\\ \bzero&\bA_{22,n}^{-1}
\end{pmatrix}
\begin{pmatrix}
\bzero\\\bar  \bY_1-\bar \bY_2 
\end{pmatrix}.
\end{align*}
It may be shown using (\ref{SigmahatinvA}) and (\ref{ASigmahatinvA}) that 
$$
\begin{pmatrix}\bzero &\bzero\\ \bzero&\bA_{22,n}^{-1}\end{pmatrix}\bA\bSigmahat^{-1}\bA\bSigmahat^{-1}\bA\begin{pmatrix}\bzero &\bzero\\ \bzero&\bA_{22,n}^{-1}\end{pmatrix} = \begin{pmatrix}
\boldsymbol{*} & \boldsymbol{*} \\ \boldsymbol{*} & 
N^2\bA_{22,N}^{-1}\bA_{22,n}\bA_{22,N}^{-1}\end{pmatrix},
$$
where $``{\boldsymbol *}''$ denotes terms whose explicit expressions are not needed.  Consequently, 
\begin{align*}
\frac1{4\bartau^2}T_3 &= \begin{pmatrix}\bzero\\\bar  \bY_1-\bar \bY_2 
\end{pmatrix}'
\begin{pmatrix}
\boldsymbol{*}&\boldsymbol{*}\\ \boldsymbol{*}&N^2\bA_{22,N}^{-1}\bA_{22,n}\bA_{22,N}^{-1}
\end{pmatrix}
\begin{pmatrix}
\bzero &\bzero\\ \bzero&\bA_{22,n}^{-1}
\end{pmatrix}
\begin{pmatrix}
\bzero\\\bar  \bY_1-\bar \bY_2 
\end{pmatrix}\\
&=N^2(\bar\bY_1-\bar\bY_2)'\bA_{22,N}^{-1}\bA_{22,n}\bA_{22,N}^{-1}(\bar  \bY_1-\bar\bY_2) \\
&=O_p(1).
\end{align*}

As for $T_4$, we have 
\begin{align*}
\frac1{\bartau^2}T_4 & = 
\sum_{j=1}^n (\bZ_j-\bar\bZ_1)'\bSigmahat^{-1}(\bZ_j-\bar\bZ_1){\wtilde\bY}'\bSigmahat^{-1}\wtilde\bY \\
& = {\wtilde\bY}'\bSigmahat^{-1}\wtilde\bY \cdot \tr\,\bSigmahat^{-1}\sum_{j=1}^n(\bZ_j-\bar\bZ_1)(\bZ_j-\bar\bZ_1)'\\
& = N(\bar\bY_1-\bar\bY_2)'\bA_{22,N}^{-1}(\bar\bY_1-\bar\bY_2)\, \tr\,\bSigmahat^{-1}\bA.
\end{align*}
Applying (\ref{SigmahatinvA}), we obtain 
$$
\frac1{\bartau^2}T_4 = N(\bar\bY_1-\bar\bY_2)'\bA_{22,N}^{-1}(\bar\bY_1-\bar\bY_2)(n\, \tr\bI_p+N\, \tr\bA_{22,N}^{-1}\bA_{22,n}) = O_p(1).
$$
In the case of $T_5$, we find that  
$$
T_5 = n\bartau^2({\wtilde\bY}'\bSigmahat^{-1}\wtilde\bY)^2 = nN\bartau^2[(\bar\bY_1-\bar\bY_2)'\bA_{22,N}^{-1}(\bar\bY_1-\bar\bY_2)]^2 = O_p(N^{-2})
$$
because $\bA_{22,N}^{-1} = O_p(N^{-1})$ and $(\bar\bY_1-\bar\bY_2)(\bar\bY_1-\bar\bY_2)' = O_p(N^{-1})$.

Collecting together these five expansions, we deduce that 
\begin{align}
\frac1N b_{2,p,q}^{(1)} = & N^{-1}\sum_{j=1}^5 T_j \nonumber \\
= & N^{-1} \Big\{\sum_{j=1}^n[(\bZ_j-\bar\bZ_1)'\bSigmahat^{-1}(\bZ_j-\bar\bZ_1)]^2\nonumber\\
&+4\bartau\sum_{j=1}^n(\bZ_j-\bar\bZ_1)'\bSigmahat^{-1}(\bZ_j-\bar\bZ_1)(\bZ_j-\bar\bZ_1)'\bSigmahat^{-1}\wtilde\bY\nonumber\\
&+4\bartau^2N^2(\bar  \bY_1-\bar \bY_2)'\bA_{22,N}^{-1}\bA_{22,n}\bA_{22,N}^{-1}(\bar  \bY_1-\bar \bY_2)
\nonumber\\&+\bartau^2N(\bar\bY_1-\bar\bY_2)'\bA_{22,N}^{-1}(\bar\bY_1-\bar\bY_2)(np+N\tr\bA_{22,N}^{-1}\bA_{22,n})
\nonumber\\
&+\bartau^3nN^2[(\bar\bY_1-\bar\bY_2)'\bA_{22,N}^{-1}(\bar\bY_1-\bar\bY_2)]^2\Big\}\nonumber\\
= & N^{-1}\sum_{j=1}^n(\bZ_j'\bZ_j)^2-2N^{-1/2}\tau \, \tr(\wtilde\bB\bXi)-4\tau\bTheta'\bar\bZ + 4\tau\bartau\bTheta'\wtilde\bY^* + o(N^{-1/2}).
\label{asymptT1}
\end{align}
Finally, it follows from (\ref{asymptT2}) and (\ref{asymptT1}) that 
\begin{align}
\label{asymptT}
b_{2,p,q} = \, & N^{-1} \(c_1 b_{2,p,q}^{(1)} + c_2 b_{2,p,q}^{(2)}\) \nonumber \\
= \, & c_1\Big\{N^{-1}\sum_{j=1}^n(\bZ_j '\bZ_j )^2-2N^{-1/2}\tau \, \tr(\wtilde\bB\Xi )-4\tau\bTheta'\bar\bZ_1 + 4\tau\bartau \bTheta'\wtilde\bY^*\Big\} \nonumber \\ 
& + c_2\Big\{N^{-1} \sum_{j=n+1}^N (\bY_j'\bY_j)^2-2N^{-1/2}\bartau \, \tr(\wtilde\bB_{22}\Xi^*)-4\bartau{\bTheta^*}'\bar\bY\Big\} + o(N^{-1/2}).
\end{align}

\section{The asymptotic distribution of \texorpdfstring{$\boldsymbol{b_{2,p,q}}$}{b}}
\label{basymptotics}
\setcounter{equation}{0}

\subsection{The proof of Theorem \ref{basympthm}}

To obtain the asymptotic distribution of $b_{2,p,q}$, we first simplify the expression (\ref{asymptT}).  Recall that 
$$
\bB^* = \begin{pmatrix}\bB_{11}^*&\bB_{12}^*\\\bB_{21}^*&\bB_{22}^*\end{pmatrix}
$$
and 
$$
\wtilde\bB = \begin{pmatrix}\bB_{11}^*&\bB_{12}^*\\\bB_{21}^*&\wtilde\bB_{22}\end{pmatrix}
$$ 
where
\begin{align*}
\begin{array}{ll}
\bB_{11}^*= N^{1/2} \Big(\dfrac1n\sum\limits_{j=1}^n\bX_j\bX_j'-\bI_{p}\Big),&
\bB_{12}^*=\bB_{21}'= N^{1/2} \dfrac1n\sum\limits_{j=1}^n\bX_j\bY_j',\\
\bB_{22}^*= N^{1/2} \Big(\dfrac1n\sum\limits_{j=1}^n\bY_j\bY_j'-\bI_{q}\Big),&
\wtilde\bB_{22}= N^{1/2} \Big(\dfrac1N\sum\limits_{j=1}^N\bY_j\bY_j'-\bI_{q}\Big).
\end{array}
\end{align*}
Define 
$$
\wtilde\bB_{22}^* = \sqrt N \Big(\frac1{N-n}\sum_{j=n+1}^N\bY_j\bY_j'-\bI_{q}\Big),
$$
then $\wtilde\bB_{22}=\tau\bB_{22}^*+ \bartau\wtilde\bB_{22}^*$.  Recalling from (\ref{Xidef}) and (\ref{Xi}) the definitions of the matrix $\bXi$ and its components $\bXi_{ij}$, direct algebraic calculations reveal that 
\begin{align*}
c_1\tau \, \tr(\wtilde\bB\bXi )&+c_2\bartau \, \tr(\wtilde\bB_{22}\bXi^*)\\
= \ & c_1\tau\{\, \tr(\bB_{11}^*\bXi_{11}+\bB_{12}^*\bXi_{21})+\, \tr(\bB_{21}^*\bXi_{12}+(\tau\bB_{22}^*+ \bartau\wtilde\bB_{22}^*)\bXi_{22})\}\\
& + c_2\bartau \, \tr((\tau\bB_{22}^*+ \bartau\wtilde\bB_{22}^*)\Xi^*)\\
= \ & c_1\tau\{\, \tr(\bB_{11}^*\bXi_{11}+\bB_{12}^*\bXi_{21})+\, \tr(\bB_{21}^*\bXi_{12}+\tau\bB_{22}^*\bXi_{22})\}+c_2\bartau \, \tr(\tau\bB_{22}^*\bXi^*)\\
& + c_1\tau \, \tr( \bartau\wtilde\bB_{22}^*\bXi_{22})+c_2\bartau \, \tr(\bartau\wtilde\bB_{22}^*\Xi^*)\\
\equiv \ &\tau \tr\bB^*\wtilde \bXi+\bartau \tr \wtilde\bB_{22}^*\wtilde\bXi_{22},
\end{align*}
Since $\bar\bZ_1 = \begin{pmatrix}\bar\bX \\ \bar\bY_1\end{pmatrix}$, $\wtilde\bY^* = \begin{pmatrix}\bzero \\ \bar\bY_1-\bar\bY_2\end{pmatrix}$, and $\bar\bY = \tau\bar\bY_1+\bartau\bar\bY_2$, then it follows from the definitions of $\bTheta$ and $\wtilde\bTheta$ in (\ref{Theta}) and (\ref{Thetatilde}), respectively, that 
\begin{align*}
c_1\tau\bTheta'\bar\bZ_1 - c_1\tau\bartau \bTheta'\wtilde\bY^* & + c_2\bartau{\bTheta^*}'\bar\bY \\
= \ &
c_1\tau(\bTheta_1'\bar\bX+\bTheta_2'\bar\bY_1) -c_1\tau\bartau \bTheta_2'(\bar\bY_1-\bar\bY_2)+c_2\bartau{\bTheta^*}'(\tau\bar\bY_1+\bartau\bar\bY_2)\\
= \ &
c_1\tau\bTheta_1'\bar\bX+(c_1\tau\bTheta_2-c_1\tau\bartau \bTheta_2+c_2\tau\bartau{\bTheta^*})'\bar\bY_1+(c_1\tau\bartau \bTheta_2+c_2\bartau^2{\bTheta^*})'\bar\bY_2 \\
= \ &
c_1\tau\bTheta_1'\bar\bX+\tau(c_1\tau\bTheta_2+c_2\bartau{\bTheta^*})'\bar\bY_1+\bartau(c_1\tau \bTheta_2+c_2\bartau{\bTheta^*})'\bar\bY_2 \\
\equiv \ & \tau\wtilde\bTheta'\bar\bZ_1+\bartau\wtilde\bTheta_2'\bar\bY_2 .
\end{align*}
Hence, 
\begin{align*}
b_{2,p,q} &= \frac{c_1}{N} \sum_{j=1}^n(\bZ_j'\bZ_j)^2 + \frac{c_2}{N} \sum_{j=n+1}^N(\bY_j'\bY_j)^2 \\
& \qquad -2N^{-1/2}(\tau \tr\bB^*\wtilde \bXi+\bartau \tr \wtilde\bB_{22}^*\wtilde\bXi_{22}) - 4(\tau\wtilde\bTheta'\bar\bZ_1+\bartau\wtilde\bTheta_2'\bar\bY_2) + o(N^{-1/2}).
\end{align*}

Since $\tau = n/N$ and $\bartau = (N-n)/N$ then $N^{-1/2} = n^{-1/2}\tau^{1/2}$ and $N^{-1/2} = (N-n)^{-1/2}\bartau^{1/2}$.  Define 
\begin{align}
\label{UVW}
\begin{array}{ll}
\begin{pmatrix}U_{1}\\V_{1}\\\bW_{1}\end{pmatrix} = n^{-1/2}\sum\limits_{j=1}^n\wtilde\bZ_j, & \quad
\wtilde\bZ_j = \begin{pmatrix}\|\bZ_j\|^4-E\|\bZ\|^4\\\bZ_j'\wtilde \bXi\bZ_j-E(\bZ'\wtilde \bXi\bZ)\\\bX_{j}\\\bY_j\end{pmatrix},\\
 & \\
\begin{pmatrix}U_{2}\\V_{2}\\\bW_{2}\end{pmatrix} = (N-n)^{-1/2} \sum\limits_{j=n+1}^N\wtilde\bZ_j^*, & \quad
\wtilde\bZ_j^* = \begin{pmatrix}\|\bY_j\|^4-E\|\bY\|^4\\\bY_j'\wtilde\bXi\bY_j-E(\bY'\wtilde \bXi_{22}\bY)\\ \bY_j\end{pmatrix};
\end{array}
\end{align}
then, 
\begin{align*}
N^{1/2}\big(&b_{2,p,q}-c_1\tau E|\bZ|^4-c_2\bar\tau E|\bY|^4\big) \\
= \ & N^{-1/2}\Big\{c_1\sum_{j=1}^n (\bZ_j'\bZ_j)^2+c_2\sum_{j=n+1}^N(\bY_j'\bY_j)^2\Big\} - N^{1/2}\Big\{c_1\tau E|\bZ|^4-c_2\bar\tau E|\bY|^4\Big\} \\
& - 2N^{-1/2}(\tau \tr\bB^*\wtilde \bXi+\bartau \tr \wtilde\bB_{22}^*\wtilde\bXi_{22}) - 4(\tau\wtilde\bTheta'\bar\bZ_1+\bartau\wtilde\bTheta_2'\bar\bY_2) + o(N^{-1/2}) \\
= \ & \tau^{1/2}n^{-1/2}c_1\sum_{j=1}^n (\bZ_j'\bZ_j)^2+\bar\tau^{1/2}{(N-n)}^{-1/2}c_2\sum_{j=n+1}^N(\bY_j'\bY_j)^2 \\
& - n^{1/2}c_1\tau^{1/2} E|\bZ|^4-(N-n)^{1/2}c_2\bar\tau^{1/2} E|\bY|^4\\
& - 2N^{-1/2}(\tau \tr\bB^*\wtilde \bXi+\bartau \tr \wtilde\bB_{22}^*\wtilde\bXi_{22}) - 4(\tau\wtilde\bTheta'\bar\bZ_1+\bartau\wtilde\bTheta_2'\bar\bY_2) + o(N^{-1/2}).
\end{align*}
Writing each $(\bZ_j'\bZ_j)^2 = (\bZ_j'\bZ_j)^2 - E|\bZ|^4 + E|\bZ|^4$, and similarly for each $|\bY|^4$, we obtain 
\begin{align*}
N^{1/2}\big(&b_{2,p,q}-c_1\tau E|\bZ|^4-c_2\bar\tau E|\bY|^4\big) \\
= \ & \tau^{1/2}c_1n^{-1/2}\sum_{j=1}^n\((\bZ_j'\bZ_j)^2-E|\bZ|^4\) 
+ \bartau^{1/2}{(N-n)}^{-1/2}c_2\sum_{j=n+1}^N\((\bY_j'\bY_j)^2- E|\bY|^4\)\\
& - 2N^{-1/2}(\tau \tr\bB^*\wtilde \bXi+\bartau \tr \wtilde\bB_{22}^*\wtilde\bXi_{22}) - 4(\tau\wtilde\bTheta'\bar\bZ_1+\bartau\wtilde\bTheta_2'\bar\bY_2) + o(N^{-1/2}) \\
= \ & \tau^{1/2}c_1U_1+\bar\tau^{1/2}c_2U_2\\
& - 2N^{-1/2}(\tau \tr\bB^*\wtilde \bXi+\bartau \tr \wtilde\bB_{22}^*\wtilde\bXi_{22}) - 4(\tau\wtilde\bTheta'\bar\bZ_1+\bartau\wtilde\bTheta_2'\bar\bY_2) + o(N^{-1/2}),
\end{align*}
so we obtain 
\begin{align*}
N^{1/2}(b_{2,p,q}-c_1\tau E\|\bZ\|^4-c_2\bartau E\|\bY\|^4)
= \ & c_1\tau^{1/2}U_1-2\tau^{1/2}V_1-4 \tau^{1/2}\wtilde\bTheta'\bW_{1}\\
& + \, c_2\bartau^{1/2}U_2-2\bartau^{1/2}V_2-4 \bartau^{1/2}\wtilde\bTheta_2'\bW_{2}+o(1)\\
\equiv \ & \wtilde b_{2,p,q}^{(1)}+\wtilde b_{2,p,q}^{(2)}+o(1),
\end{align*}
where 
\begin{equation}
\label{b1b2}
\begin{aligned}
\wtilde b_{2,p,q}^{(1)} = \ & c_1\tau^{1/2}U_1-2\tau^{1/2}V_1-4 \tau^{1/2}\wtilde\bTheta'\bW_{1}, \\
\wtilde b_{2,p,q}^{(2)} = \ & c_2\bartau^{1/2}U_2-2\bartau^{1/2}V_2-4 \bartau^{1/2}\wtilde\bTheta_2'\bW_{2}.
\end{aligned}
\end{equation}
Note that $\wtilde b_{2,p,q}^{(1)}$ depends on $\bZ_1,\ldots,\bZ_n$ only, and $\wtilde b_{2,p,q}^{(2)}$ depends on $\bY_{n+1},\ldots,\bY_N$ only.  Since $\bZ_1,\ldots,\bZ_n$ and $\bY_{n+1},\ldots,\bY_N$ are independent, then $\wtilde b_{2,p,q}^{(1)}$ and $\wtilde b_{2,p,q}^{(2)}$ also are independent.  

By (\ref{UVW}), $\wtilde\bZ_1,\ldots,\wtilde\bZ_n$ are mutually independent and identically distributed with $E(\wtilde\bZ_1) = \bzero$ and covariance matrix
\begin{align}
E(\wtilde\bZ_1\wtilde\bZ_1') = 
\begin{pmatrix}
\Var(\|\bZ\|^4) & \Cov(\|\bZ\|^4,\bZ'\wtilde\bXi\bZ) & E(\|\bZ\|^4\,\bZ')\\
\Cov(\|\bZ\|^4,\bZ'\wtilde\bXi\bZ) & \Var(\bZ'\wtilde\bXi\bZ) & E(\bZ'\wtilde\bXi\bZ\bZ')\\
E(\|\bZ\|^4\,\bZ) & E(\bZ\bZ'\wtilde\bXi\bZ) & \bI_{p+q}
\end{pmatrix},
\end{align}
and 
$\wtilde\bZ_{n+1}^*$,  $\cdots$, $\wtilde\bZ_N^*$ are mutually independent and identically distributed with $E(\wtilde\bZ_N^*) = \bzero$ and covariance matrix
\begin{align}
E(\wtilde\bZ_N^* \wtilde\bZ_N^*{}') = 
\begin{pmatrix}
\Var(\|\bY\|^4) & \Cov(\|\bY\|^4,\bY'\wtilde\bXi_{22}\bY) & E(\|\bY\|^4\,\bY')\\
\Cov(\|\bY\|^4,\bY'\wtilde\bXi_{22}\bY) & \Var(\bY'\wtilde\bXi_{22}\bY) & E(\bY'\wtilde\bXi_{22}\bY\bY')\\
E(\|\bY\|^4\,\bY)&E(\bY\bY'\wtilde\bXi_{22}\bY)&\bI_{p+q}
\end{pmatrix}.
\end{align}
Since $\wtilde b_{2,p,q}^{(1)}$ and $\wtilde b_{2,p,q}^{(2)}$ are independent, it follows that 
\begin{eqnarray*}
\Var(b_{2,p,q}) & = & \Var\big(\,\wtilde b_{2,p,q}^{(1)}\big) + \Var\big(\,\wtilde b_{2,p,q}^{(2)}\big)\\ 
& = & \tau\bnu_1'E(\wtilde\bZ_1\wtilde\bZ_1')\bnu_1 + \bartau\bnu_2'E(\wtilde\bZ_N^*\wtilde\bZ_N^*{}')\bnu_2,
\end{eqnarray*}
where 
\begin{align*}
\bnu_1 = 
\begin{pmatrix}
c_1 \\ -2 \\ -4\wtilde \bTheta
\end{pmatrix},
\quad 
\bnu_2 =
\begin{pmatrix}
c_2 \\ -2 \\ -4\wtilde \bTheta_2
\end{pmatrix}.
\end{align*}
Applying the Central Limit Theorem to (\ref{UVW}) and (\ref{b1b2}), we obtain Theorem \ref{basympthm}.  $\qed$

\subsection{The proof of Corollary \ref{basympcor}}

To establish Corollary \ref{basympcor}, we need to calculate $\nu$ and $\sigma^2$ in (\ref{nu}) and (\ref{sigmasq}), respectively, for the null case in which the population is $N_{p+q}(\bzero,\bI_{p+q})$.   In this case, $\bX \sim N_p(\bzero,\bI_p)$, $\bY \sim N_q(\bzero,\bI_q)$ and $\bX$ and $\bY$ are independent.  Hence, $\|\bX\|^2 \sim \chi^2_p$, $\|\bY\|^2 \sim \chi^2_q$, and $\|\bZ\|^2 = \|\bX\|^2 + \|\bY\|^2 \sim \chi^2_{p+q}$.  Therefore, $E(\|\bY\|^4) = q(q+2)$, $E(\|\bZ\|^4) = (p+q)(p+q+2)$, and by substituting these results into (\ref{nu}), we obtain (\ref{nullnu}).  

Next, we calculate $\sigma^2$.  Since $\bZ \eqdist -\bZ$, it follows from a change of sign that $E(\|\bZ\|^4\,\bZ') = E(\bZ'\wtilde\bXi\bZ\bZ') = \bzero$.  Similarly, $E(\|\bY\|^4\bY') = E(\bY'\wtilde\bXi_{22}\bY\bY') = \bzero$.  Therefore, (\ref{sigmasq}) reduces to 
\begin{eqnarray*}
\sigma^2 & = & \frac{c_1^2}{\tau}\Var(\|\bZ\|^4)
+4\tau\Var(\bZ'\wtilde\bXi\bZ) 
+16\tau\wtilde\bTheta'\wtilde\bTheta 
-4c_1\Cov(\|\bZ\|^4, \bZ'\wtilde\bXi\bZ) \\
& & +\frac{c_2^2}{\bartau}\Var(\|\bY\|^4) 
+4\bartau\Var(\bY'\wtilde\bXi_{22}\bY) 
+16\bartau\wtilde\bTheta_2'\wtilde\bTheta_2 
-4c_2\Cov(\|\bY\|^4, \bY'\wtilde\bXi_{22}\bY).
\end{eqnarray*}
Note that 
\begin{align*}
\Var(\|\bY\|^4) & = E(\|\bY\|^8) - (E\|\bY\|^4)^2 \\
& = 8q(q+2)(q+3)
\end{align*}
and, similarly, $\Var(\|\bZ\|^4) = 8(p+q)(p+q+2)(p+q+3)$.  

Denote the eigenvalues of $\wtilde\bXi_{22}$ by $\lambda_j(\wtilde\bXi_{22})$, $j=1,\ldots,q$.  The distribution of $\bY$ being orthogonally invariant, we apply an orthogonal transformation to deduce that $\bY'\wtilde\bXi_{22}\bY \eqdist \sum_{j=1}^q \lambda_j(\wtilde\bXi_{22}) y_j^2$, where $y_1,\ldots,y_q$ are the independent $N(0,1)$-distributed components of $\bY$.  Then we obtain $\Var(\bY'\wtilde\bXi_{22}\bY) = 2 \, \tr(\wtilde\bXi_{22}^2)$ and $\Var(\bZ'\wtilde\bXi\bZ) = 2 \, \tr(\wtilde\bXi^2)$.  

As for the covariance terms, we again apply an orthogonal transformation to obtain 
\begin{eqnarray*}
\Cov(\|\bY\|^4,\bY'\wtilde\bXi_{22}\bY) & = & 
\Cov(\|\bY\|^4,\sum_{j=1}^q \lambda_j(\wtilde\bXi_{22})y_j^2) \\
& = & \sum_{j=1}^q \lambda_j(\wtilde\bXi_{22}) \Cov(\|\bY\|^4,y_j^2) \\
& = & \tr(\wtilde\bXi_{22}) \Cov(\|\bY\|^4,y_1^2),
\end{eqnarray*}
where we have also used the exchangeability of $y_1,\ldots,y_q$ to deduce that $\Cov(\|\bY\|^4,y_j^2) = \Cov(\|\bY\|^4,y_1^2)$ for all $j=1,\ldots,q$.  Since $\|\bY\|^2 \sim \chi^2_q$ then 
\begin{eqnarray*}
\Cov(\|\bY\|^4,y_1^2) & = & q^{-1} \sum_{j=1}^q \Cov(\|\bY\|^4,y_j^2) \\
& = & q^{-1} \Cov(\|\bY\|^4,\sum_{j=1}^q y_j^2) \\
& = & q^{-1} \Cov(\|\bY\|^4,\|\bY\|^2) = 4(q+2).
\end{eqnarray*}
Hence, $\Cov(\|\bY\|^4,\bY'\wtilde\bXi_{22}\bY) = 4(q+2) \, \tr(\wtilde\bXi_{22})$ and, similarly, $\Cov(\|\bZ\|^4,\bZ'\wtilde\bXi\bZ) = 4(p+q+2) \, \tr(\wtilde\bXi)$.  

Collecting together these results, we obtain 
\begin{align*}
\sigma^2 = \ & \frac{c_1^2}{\tau}8(p+q)(p+q+2)(p+q+3)
+8\tau \, \tr(\wtilde\bXi^2) 
+16\tau\wtilde\bTheta'\wtilde\bTheta 
-16c_1 (p+q+2) \, \tr(\wtilde\bXi) \\
& +\frac{c_2^2}{\bartau}8q(q+2)(q+3) 
+8\bartau \, \tr(\wtilde\bXi_{22}^2) 
+16\bartau\wtilde\bTheta_2'\wtilde\bTheta_2 
-16c_2 (q+2) \, \tr(\wtilde\bXi_{22}).
\end{align*}

Finally, we calculate $\wtilde\bXi$, $\wtilde\bTheta$, $\wtilde\bXi_{22}$, and $\wtilde\bTheta_2$.  By (\ref{Xi}), 
\begin{eqnarray}
\label{Xinormal}
\begin{pmatrix}
\bXi_{11}&\bXi_{12}\\\bXi_{21}&\bXi_{22}
\end{pmatrix} 
& = & 
\begin{pmatrix}
E(\|\bZ\|^2\,\bX\bX')&E(\|\bZ\|^2\,\bX\bY') \\
E(\|\bZ\|^2\,\bY\bX')&E(\|\bZ\|^2\,\bY\bY')
\end{pmatrix} \nonumber \\
& = & 
\begin{pmatrix}
E(\|\bZ\|^2\,\bX\bX')&\bzero\\\bzero&E(\|\bZ\|^2\,\bY\bY')
\end{pmatrix}.
\end{eqnarray}
By (\ref{Xitilde}), 
\begin{equation}
\label{Xitildenormal}
\wtilde\bXi = 
\begin{pmatrix}
c_1\bXi_{11}&c_1\bXi_{12}\\c_1\bXi_{21}&c_1\tau\bXi_{22}+c_2\bartau\bXi^*
\end{pmatrix} = 
\begin{pmatrix}
c_1\bXi_{11}&\bzero\\\bzero&c_1\tau\bXi_{22} + c_2\bartau\bXi^*
\end{pmatrix}.
\end{equation}
By (\ref{Xistar}), 
\begin{eqnarray*}
\bXi^* = E((\bY\bY')^2) 
& = & E(\bY(\bY'\bY)\bY') \\
& = & E(\|\bY\|^2\bY\bY') = \big(E(\|\bY\|^2 Y_iY_j)\big).
\end{eqnarray*}
If $i \neq j$ then, by the exchangeability of $Y_1,\ldots,Y_q$, we have 
\begin{eqnarray*}
E(\|\bY\|^2 Y_iY_j) & = & E(\|\bY\|^2 Y_1Y_2) \\
& = & E((Y_1^2+Y_2^2+Y_3^2+\cdots+Y_q^2) Y_1Y_2) = 0,
\end{eqnarray*}
the last equality following from the mutual independence of $Y_1,\ldots,Y_q$ and the fact that each has mean $0$.  If $i=j$ then by exchangeability, 
\begin{eqnarray*}
E(\|\bY\|^2 Y_iY_j) & = & E(\|\bY\|^2 Y_j^2) \\
& = & \frac1q E(\|\bY\|^2 \cdot \sum_{j=1}^q Y_j^2) = \frac1q E(\|\bY\|^4) = q+2;
\end{eqnarray*}
therefore, $\bXi^* = (q+2)\bI_q$.  Next, by (\ref{Xinormal}), 
\begin{eqnarray*}
\bXi_{22} = E(\|\bZ\|^2\,\bY\bY') 
& = & E((\|\bX\|^2+\|\bY\|^2)\,\bY\bY') \\
& = & E(\|\bX\|^2) E(\bY\bY') + E(\|\bY\|^2\,\bY\bY') 
= (p+q+2)\bI_q,
\end{eqnarray*}
and by interchanging the roles of $\bX$ and $\bY$ in this latter calculation, we obtain $\bXi_{11} = E(\|\bZ\|^2\,\bX\bX') = (p+q+2)\bI_p$.  Inserting these results at (\ref{Xitildenormal}), we obtain 
$$
\wtilde\bXi = 
\begin{pmatrix} 
c_1(p+q+2)\bI_p & \bzero \\ \bzero & (c_1\tau(p+q+2)+c_2\bartau(q+2))\bI_q 
\end{pmatrix}
$$
so that $\wtilde\bXi_{22} = (c_1\tau(p+q+2)+c_2\bartau(q+2))\bI_q$.  

As for $\wtilde \bTheta$ and $\wtilde\bTheta_2$, it follows from (\ref{Theta}) and (\ref{Thetatilde}) that 
$$
\wtilde\bTheta \equiv \begin{pmatrix}\wtilde\bTheta_1 \\ \wtilde\bTheta_2 \end{pmatrix} = \begin{pmatrix}c_1\bTheta_1 \\ c_1\tau\bTheta_2+c_2\bartau\bTheta^*\end{pmatrix}
$$
where 
$$
\bTheta = \begin{pmatrix}\bTheta_1 \\ \bTheta_2 \end{pmatrix} = \begin{pmatrix} E(\|\bZ\|^2\bX) \\ E(\|\bZ\|^2\bY)\end{pmatrix} = \bzero.
$$
By (\ref{Xistar}), $\bTheta^* = E(\|\bY\|^2 \bY) = \bzero$; hence, $\wtilde\bTheta = \bzero$.  Inserting these results in (\ref{sigmasq}), we obtain (\ref{nullsigmasq}).  
$\qed$

\end{document}